\documentclass[10pt,a4paper,leqno]{amsart}
\usepackage[english]{babel}
\usepackage[T1]{fontenc}
\usepackage{amssymb}
\usepackage{amsfonts}
\usepackage{amsmath}
\usepackage{mathptm}

\usepackage{eucal}
\begin{document}
\thispagestyle{empty}
\date{September 16, 2010}
\title{{{On an inequality related to the radial growth of quasinearly subharmonic functions in locally uniformly homogeneous spaces}}}
\author{Juhani Riihentaus}
\address{\begin{tabular}{cccccc}
Department of Mathematical Sciences&&&&Department of Physics and Mathematics\\
University of Oulu&&&&University of Eastern Finland\\
P.O. Box 3000&&&&P.O. Box 111\\
FI-90014 University of Oulu, Finland&&&&FI-80101 Joensuu, Finland\\
\end{tabular}}
\email{juhani.riihentaus@gmail.com \, {\textrm{ and }} \, juhani.riihentaus@uef.fi}
\begin{abstract} We begin by recalling the definition of  nonnegative   quasinearly subharmonic functions on locally uniformly homogeneous spaces. Recall that these spaces and this function class are rather general:   Among others subharmonic, quasisubharmonic and  nearly subharmonic functions  on domains of  Euclidean spaces ${\mathbb{R}}^n$, $n\geq 2$, are included.  
The following result of  Gehring and Hallenbeck is classical: Every
subharmonic function, defined  and ${\mathcal{ L}}^p$-integrable for some $p$,
$0<p<+\infty$,
on the unit disk ${\mathbb{D}}$ of the complex plane ${\mathbb{C}}$ is for almost all $\theta$ of the form $o((1-\vert z\vert
)^{-1/p})$, uniformly as
$z\to e^{i\theta}$ in any Stolz domain. 
Recently both Pavlovi\'c and Riihentaus have given related and partly more general results on domains of ${\mathbb{R}}^n$, $n\geq 2$. Now we extend one of these  results  to quasinearly subharmonic functions on locally uniformly homogeneous spaces.
\end{abstract} 
\maketitle
\section{Introduction}
\subsection{Locally uniformly homogeneous spaces} The definition of locally uniformly homogeneous spaces was given in \cite{PaRi10}. However, for the convenience of the reader we recall it  here, too. A set $X$ is \emph{a locally uniformly homogeneous space} if the following conditions  are satisfied:
\begin{itemize}
\item[(i)] $X$ is a topological space.
\item[(ii)] There is a Borel measure $\mu $ defined on $X$.
\item[(iii)] There is a \emph{quasimetric} (\emph{quasidistance}) on $X$, that is, there is a constant $K\geq 1$ and a mapping $d_K: X\times X\rightarrow [0,+\infty )$ such that 
\begin{itemize}
\item[$1^o$] $d_K(x,y)=d_K(y,x)$ for all $x,y\in X$,
\item[$2^o$] $d_K(x,y)=0$ if and only if $x=y$,
\item[$3^o$] $d_K(x,y)\leq K[d_K(x,z)+d_K(z,y)]$ for all $x,y,z\in X$,
\item[$4^o$] the (($K$-)quasi)balls $B_K(x,r)$,
\[B_K(x,r):=\{\,y\in X\, :\, d_K(x,y)<r\,\},\]
centered at $x$ and of radii $r>0$, form a basis of open neighborhoods at the point $x\in X$,
\item[$5^o$] $0<\mu (B_K(x,r))<+\infty$ for all $x\in X$ and $r>0$,
\item[$6^o$] there exist absolute constants $A=A(K)\geq 1$ and $\rho _0=\rho_0(K)>0$ such that 
\[\mu (B_K(x,r))\leq A\, \mu (B_K(x,\frac{r}{2}))\]
for all $x\in X$ and all $r$, $0<r\leq \rho _0$.
\end{itemize}
\end{itemize}
\subsubsection{Remark} Locally uniformly homogeneous spaces are slightly more general than  spaces of homogeneous type,  defined and considered by Coifman and Weiss  \cite{CoWe71}, pp.~66--68, and \cite{CoWe77}, pp.~587--590. As a matter of fact, the only difference with their definition is that, instead of  the  above \mbox{condition $6^o$,} Coifman and Weiss use  the  stronger condition:
\begin{itemize}
\item[$6'^o$] There exists an absolute  constant $A=A(K)\geq 1$ such that 
\[\mu (B_K(x,r))\leq A\, \mu (B_K(x,\frac{r}{2}))\]
for all $x\in X$ and all  $r>0$.
\end{itemize}
For a  list of examples of spaces of homogeneous type, see \cite{CoWe77}, pp.~588--590.
\subsubsection{Remark}  In order to be able to consider Hausdorff measures on locally uniformly homogeneous spaces, we make the following additional assumption (cf. \cite{Ma95}, p.~54): Let $X$ be a locally uniformly homogeneous space. Suppose that $X$ satisfies  the following additional condition: 
\begin{equation}
\emph{For every $\delta >0$ there are $E_j\in {\mathcal{F}}$, $j=1,2, \dots$, such that $X=\cup_{j=1}^{+\infty}E_j$ and $d_K(E_j)\leq \delta$,} 
\end{equation} 
\emph{where ${\mathcal{F}}=\{\, B_K(x,r)\,:\, x\in X,\, r>0\,\}$.} 
Then for each $d>0$ one can define   in $X$ a $d$-dimensional Hausdorff (outer) measure ${\mathcal{H}}_K^d$,  which is  a ($K$-quasi)metric (outer) measure in the following sense: If $A, B\subset X$ such that $d_K(A,B)>0$, then ${\mathcal{H}}_K^d(A\cup B)={\mathcal{H}}_K^d(A)+{\mathcal{H}}_K^d(B)$. As a matter of fact, in the standard definition (see e.g.   \cite{LuMa95}, pp.~125--126), just work with the quasimetric $d_K$ instead of the  metric  $d$ (or $\rho$). One sees also, that all Borel sets of $X$ are ${\mathcal{H}}_K^d$-measurable. Above we have used the following notation: If $A,B\subset X$, then  \[d_K(A):=\sup \{\, d_K(x,y)\, :\, x,y\in A\, \} \quad {\textrm{and}}\quad  d_K(A,B):=\inf \{\, d_K(x,y)\, :\, x\in A,\, y\in B\, \}.\]
\subsection{Quasinearly subharmonic functions} Though the definition of quasinearly subharmonic functions in locally uniformly homogeneous spaces was given in \cite{PaRi10},  we recall it also here for the convenience of the reader. Let $X$ be a locally uniformly homogeneous space.  Let $u:\,X\rightarrow [0,+\infty )$ be Borel measurable. Let $C\geq 1$. Then $u$ is $C$\emph{-quasinearly subharmonic in}  $X$ if there is a constant $\epsilon _0=\varepsilon_0 (u)$ (depending on the considered function $u$), $0<\varepsilon_0<1$, such that 
for each open set $\Omega \subset X$, $\Omega\ne X$, for each $x\in \Omega $ and each $r$, $0<r\leq \min \{\,\rho_0, \, \varepsilon_0 \delta _K^{\Omega }(x)\,\}$, one has $u\in {\mathcal{L}}^1(B_K(x,r))$ and 
\[u(x)\leq \frac{C}{\mu (B_K(x,r))}\int\limits_{B_K(x,r)}u(y)\, d\mu (y).\]
The function $u$ is \emph{quasinearly subharmonic in}  $X$ if $u$ is $C$-quasinearly subharmonic for some $C\geq 1$.  

Above (and below) we have used  the following notation: $\delta_K^{\Omega}(x)$, or shortly $\delta_K(x)$, is 
the ($K$-)quasidistance from $x\in \Omega$ to $\partial \Omega$, and thus defined by
\begin{equation*} \delta_K(x):=\delta_K^{\Omega}(x):=\inf \{\, d_K(x,y)\,:\, y\in \Omega^c \,\}  
\end{equation*}
where $\Omega^c$ is the complement of $\Omega$, taken in $X$.
\subsubsection{Examples} Quasinearly subharmonic functions, especially nearly subharmonic, quasisubharmonic and subharmonic functions in an open subset $D$ of an Euclidean space ${\mathbb{R}}^n$, $n\geq 2$, give  examples of quasinearly subharmonic functions in a locally uniformly homogeneous space. As an additional example, we recall that  $B^{2n}$, the unit ball of ${\mathbb{C}}^n$, $n\geq 1$, is locally uniformly homogeneous, and nonnegative ${\mathcal{M}}$-subharmonic functions on $B^{2n}$ (see e.g. \cite{St94}, p.~31, and \cite{St97}, p.~3774) are $1$-quasinearly subharmonic. For further examples, see \cite{PaRi10}. 

For the definition, examples  and properties of quasinearly subharmonic functions (sometimes, however, perhaps with a different terminology) in domains of an Euclidean space ${\mathbb{R}}^n$, $n\geq 2$, see e.g. \cite{Pa94}, pp.~18--19, \cite{Pa96}, pp.~15--16, \cite{Ri00}, p.~233,  \cite{Ri03}, p.~171, \cite{Ri04}, pp.~196--197, \cite{Ri06$_1$}, p.~28, \cite{Ri06$_2$}, p.~158, \cite{Ko07}, pp. 243--244, \cite{Ri07}, p.~52, \cite{PaRi08}, pp.~90--91, 
\cite{Ri08$_1$}, pp.~2--3, 
\cite{Ri08$_2$}, p.~2614, \cite{Ri09}, pp.~129--130, \cite{DoRi10$_1$}, pp.~2--6, \cite{DoRi10$_2$}, and the references therein. In this connection, see also \cite{Vu82}, pp.~259, 263.
\subsection{Weighted boundary behavior} The following theorem is a special case of  the original  result of Gehring  \cite{Ge57}, Theorem~1, p.  77, and 
of Hallenbeck \cite{Ha92}, Theorems~1 and 2, pp.~117--118,  and of the later and more general results of Stoll \cite{St98},  Theorems~1 and 2, pp.~301--302, 307:
\subsection{Theorem} {\emph{If $u$ is a function harmonic in the unit disk ${\mathbb{D}}$ of the complex plane ${\mathbb{C}}$ such that 
\begin{equation} I(u):=\int_{{\mathbb{D}}}\mid u(z)\mid^p (1-\mid z\mid )^{\beta} \, dm_2(z)<+\infty ,\end{equation}
where $p>0$, $\beta >-1$, then 
\begin{equation} \lim_{r\rightarrow 1-}\mid u(re^{i\theta })\mid ^p(1-r)^{\beta +1}=0\end{equation}
for almost all $\theta \in [0,2\pi )$. Above $m_2$ is the Lebesgue measure in ${\mathbb{R}}^2$.}}

Observe that  Gehring, Hallenbeck and Stoll  considered in fact subharmonic functions and  that the limit in (3) was  uniform in Stolz approach regions, in Stoll's result in even more general regions. For  more general results, see \cite{Ri00}, Theorem, p.~233, \cite{Mi01}, Theorem~2, p.~73, \cite{Ri03}, Theorem~2, pp.~175--176, \cite{Ri04}, Theorem~3.4.1, pp.~198--199, \cite{Ri06$_1$}, Theorem, p.~31, and \cite{PaRi08}, Theorem~4, p.~102.

Gehring's  proof was based on Hardy-Littlewood inequality, whereas the other authors based their proofs, more or less, on  certain generalized mean value inequalities for subharmonic functions.  
For such inequalities  and related properties, see \cite{FeSt72},  Lemma~2, p.~172, \cite{Ku74}, Theorem~1, p.~529, \cite{To86}, pp.~188-190,  \cite{Pa88}, pp.~53, 64--65, \cite{Ri89}, Lemma, p.~69,  \cite{Ha92}, Lemma~1, p.~113,  \cite{Mi01}, p.~68, \cite{Ri01}, Theorem, p.~188, and the references therein.

With the aid of the following Theorem~1.5, see \cite{Pa08}, Theorem~1, pp.~433--434,  Pavlovi\'c showed that the convergence in (3)  is dominated. 
At the same time he pointed out that  whole Theorem~1.4 follows from his result:
\subsection{Theorem} {\emph{If $u$ is a function harmonic  in ${\mathbb{D}}$ satisfying}} (2){\emph{, where $p>0$, $\beta >-1$, then 
\begin{equation*} J(u):=\int\limits_{0}^{2\pi}\sup_{0<r<1}\mid u(re^{i\theta })\mid ^p (1-r)^{\beta +1} \, d\theta <+\infty .\end{equation*}
Moreover, there is a constant $C=C_{p,\beta }$ such that $J(u)\leq C\,I(u)$.}}

In \cite{Ri09}, Theorems~1 and 2, pp.~131--132,   we extended Theorem~1.5 to the case, where, instead of absolute values of harmonic functions in the unit disk ${\mathbb{D}}$ of  the complex plane ${\mathbb{C}}$, one considers more generally  nonnegative quasinearly subharmonic functions in rather general domains of ${\mathbb{R}}^n$, $n\geq 2$. Now our aim is to  extend this cited Theorem~1  even further: We will give a related  result for  quasinearly subharmonic functions in locally uniformly homogeneous spaces, satisfying the above additional assumption (1), see Theorem~2.5 below. As an application, we get in  Corollary~2.6 below  a weighted boundary behavior result in our rather general setup of locally uniformly homogeneous spaces.   
\subsection{Notation} Our notation is rather standard, see e.g. \cite{He71}, \cite{PaRi08} and \cite{PaRi10}. However, for the convenience of the reader 
we recall the following.
The common convention $0\cdot \infty =0$ is used. Below $X$ is always a locally uniformly homogeneous space, and  
$\Omega$ always a domain in $X$,  $\Omega \ne X$, whose boundary 
$\partial \Omega $ is Ahlfors-regular from above, with dimension $d>0$ and with constant $C_4>0$ (for the definition of this see {1.9} below).
For $\rho >0$ write $\Omega_\rho =\{\,x\in \Omega :\, \delta_K (x)<\rho \,\}.$
 $B_K(x,r)$ is the (($K$-)quasi)ball in $X$, with center $x$ and radius $r$, and 
$B_K(x)=B_K(x,\frac{1}{3K}\delta_K(x))$. 
The \mbox{$d$-dimensional} Hausdorff (outer) measure in $X$, $d>0$, constructed with the aid of the $K$-quasimetric $d_K$,  is denoted by ${\mathcal{H}}_K^d$, see Remark~1.1.2 above. 
$C_0$ and $r_0$ are fixed 
constants which  are involved with the used (and thus   fixed)  admissible function $\varphi $ (see (4) below in {1.7}). Similarly, if $\alpha >0$ is given,
$C'_1=C'_1(\alpha , C_0,K)$, $C'_2=C'_2(\alpha ,C_0, K)$ and  $C'_3=C'_3(\alpha ,C_0)$ are fixed constants, coming directly from Lemmas~2.1, 2.3 and 2.4 below. (Compare these  with the related constants $C_1$, $C_2$ and $C_3$ in \cite{Ri00}, p. 234, \cite{Ri06$_1$}, pp.~32--33,  and \cite{Ri09}, p.~129.)  
\subsection{Admissible functions}
A function $\varphi : [0,+\infty )\to [0,+\infty )$ is {\emph{admissible}}, if
it is strictly increasing, surjective, and  there are constants
$C_0=C_0(\varphi )\geq 1$ and
$r_0=r_0(\varphi )>0$ such that
\begin{equation} \varphi
(2t)\leq C_0\, \varphi
(t)
\, \, \, \,
{\textrm{ {and}}}\, \, \, \, \, \,  \varphi ^{-1} (2s)\leq C_0\, \varphi ^{-1}(s)\end{equation}
for all
$s$, $t$,  $0\leq s,\, t\leq r_0$.

Examples of admissible functions are: Functions $\varphi_1(t)=t^p$, $p>0$, nonnegative,  increasing surjective functions $\varphi_2(t)$  satisfying  the $\varDelta_2$-condition and for which the functions
$t\mapsto \frac{\varphi_2(t)}{t}$ are increasing, and functions $\varphi_3(t)=c\, t^{\alpha}[\log (\delta +t^{\gamma})]^{\beta},$
where $c>0$, $\alpha >0$, $\delta \geq 1$,  and $\beta ,\gamma \in {\mathbb {R}}$  are such that 
$\alpha +\beta \gamma >0$.
\subsection{Approach sets} Let $\varphi :[0,+\infty )\rightarrow [0,+\infty )$ be an
admissible function and let $\alpha >0$. Let $X$ be a locally uniformly homogeneous space. Let $\Omega$ be a domain in a component $X_1$ of  $X$, $\Omega \ne X_1$. For $\zeta \in \partial \Omega$   write
\[ \Gamma_{\varphi} (\zeta ,\alpha ):=\{\, x\in \Omega \, :\, \varphi (d_K(x,\zeta ))<\alpha \, \delta_K (x)\, \},\]
and  call it  a  $(\varphi ,\alpha )$-{\textit {approach set (region)}}, shortly an \emph{approach set (region)}, \emph{in}
$\Omega$ {\emph{at}} $\zeta$.  Observe that though $\partial \Omega$ is  surely nonempty, the approach set $\Gamma_{\varphi} (\zeta ,\alpha )$ may, in certain cases, be empty. Anyway, in the case of the unit disk ${\mathbb{D}}$ of the complex plane ${\mathbb{C}}$,  the choice $\varphi (t)=t$ gives  
the familiar Stolz approach regions. Choosing $\varphi (t)=t^\tau $, $\tau \geq 1$, say, one gets more general approach regions, see \cite{St98}, p.~301. 

For $x\in \Omega$ and $\alpha >0$, we  also write 
\[ \tilde{\Gamma}_{\varphi} (x,\alpha ):=\{\, \xi \in \partial \Omega \, :\, x\in \Gamma_{\varphi }(\xi ,\alpha )\, \}.\]
Moreover, for $\rho >0$ we write  
\[ {\Gamma}_{\varphi ,\rho} (\zeta ,\alpha ):=\{\, x \in \Gamma_{\varphi }(\zeta ,\alpha )\,:\, \delta_K(x)< \rho \,\}.\]

One says that $\zeta \in \partial \Omega$ is $(\varphi ,\alpha )$-\emph{accessible}, shortly \emph{accessible}, if $\Gamma_{\varphi}(\zeta ,\alpha )\cap B_K(\zeta ,\rho )\ne \emptyset$ for all $\rho >0$.
\subsection{Ahlfors-regular sets}  Let $d>0$. A set $E\subset X$ is \emph{Ahlfors-regular from above, with dimension $d$ and with  constant $C_4>0$}, shortly \emph{Ahlfors-regular from above},  if it is closed and  
\[ {\mathcal{H}}_K^d(E\cap B_K(x,r))\leq C_4r^d\]
for all $x\in E$ and $r>0$. The smallest constant $C_4$ is called the \emph{regularity constant} for $E$. 
A set $E\subset X$ is \emph{Ahlfors-regular, with dimension $d$ and with  constant $C_4>0$}, shortly \emph{Ahlfors-regular}, if it is closed  and
\[ C_4^{-1}r^d\leq {\mathcal{H}}_K^d(E\cap B_K(x,r))\leq C_4r^d\]
for all $x\in E$ and $r>0$.  

Simple examples of Ahlfors-regular sets in ${\mathbb{R}}^n$, $n\geq 2$,
are $d$-planes and $d$-dimensional Lipschitz graphs. Also certain Cantor sets and self-similar sets are Ahlfors-regular. For more details, 
see \cite{DaSe93}, pp.~9--10.   
\section{Boundary  integral inequalities}
\noindent We begin with four lemmas. Recall that  $X$ is always a  locally uniformly homogeneous space. $X_1$ will be an arbitrary  component of $X$,  and  $\Omega$  a domain in $X_1$, $\Omega \ne X_1$. Moreover $\varphi :[0,+\infty )\rightarrow [0,+\infty )$ is an admissible function, with constants $r_0$ and $C_0$. 

Let $x\in \Omega$. It is easy to see that 
\begin{equation*}
\frac{2}{3K}\delta_K(x)\leq \delta_K(y)\leq (K+\frac{1}{3})\delta_K(x)
\end{equation*}
for all $y\in B_K(x)$. 
Let $\alpha >0$. Write 
\begin{equation*}\hat{\rho}_0:= \min \{\,\frac{r_0}{2^{K+\frac{1}{3}}}, \frac{r_0}{2^{\alpha +1}},
\frac{r_0}{2^{\frac{3\alpha K}{2}C_0^{K+\frac{1}{3}}+1}(K+\frac{1}{3})},
 \frac{1}{2^{\alpha +1}}\varphi \left(\frac{r_0}{2^{K+\frac{1}{3}}}\right), {\rho}_0 \,\}.
\end{equation*} 
\subsection{Lemma} \emph{Let $\zeta \in \partial \Omega$ and $x_0\in \Gamma_{\varphi ,\rho}(\zeta , \alpha )$. Let $C_1\geq 1$ be arbitrary. Then for $C_2'=\frac{C_0^{\alpha}}{3}+KC_0^{K+\frac{1}{3}}$ and for  all $x\in B_K(x_0)$,} 
\begin{equation*} B_K(x,C_1\varphi^{-1}(\delta_K(x)))\subset B_K(x_0,C_1C_2'\varphi^{-1}(\delta_K(x_0))),\end{equation*}
\emph{provided $0<\rho \leq \hat{\rho}_0$.}
\begin{proof}
Take $z\in B_K(x,C_1\varphi^{-1}(\delta_K(x)))$ arbitrarily. Then 
\begin{align*}
d_K(x_0,z)&\leq K[d_K(x_0,x)+d_K(x,z)]<K[\frac{\delta_K(x_0)}{3K}+C_1\varphi^{-1}(\delta_K(x))]\\
&<K[\frac{\delta_K(x_0)}{3K}+C_1\varphi^{-1}((K+\frac{1}{3})\delta_K(x_0))]
<K[\frac{\delta_K(x_0)}{3K}+C_1\varphi^{-1}(2^{K+\frac{1}{3}}\delta_K(x_0))]\\
&<K[\frac{\delta_K(x_0)}{3K}+C_1C_0^{K+\frac{1}{3}}\varphi^{-1}(\delta_K(x_0))]
<K\left(\frac{C_0^{\alpha}}{3K}+C_1C_0^{K+\frac{1}{3}}\right)\varphi^{-1}(\delta_K(x_0))\leq C_1C_2'\varphi^{-1}(\delta_K(x_0)),
\end{align*}
where $C_2'=\frac{C_0^{\alpha}}{3}+KC_0^{K+\frac{1}{3}}$. Hence $z\in B_K(x_0,C_1C_2'\varphi^{-1}(\delta_K(x_0)))$. Above we have used the facts that $2^{K+\frac{1}{3}}\hat{\rho}_0\leq r_0$ and $2^{\alpha +1}\hat{\rho}_0\leq r_0$, which follow from the definition of $\hat{\rho}_0$.
\end{proof}
\subsection{Lemma} \emph{Let $\zeta \in \partial \Omega$ and $x_0\in \Gamma_{\varphi ,\rho}(\zeta , \alpha )$. Then $B_K(x_0)\subset \Gamma_{\varphi ,\rho'}(\zeta , \alpha' )$, where $\rho'=(K+\frac{1}{3})\rho$ and $\alpha'=\frac{3\alpha K}{2}C_0^{K+\frac{1}{3}}$, provided $0<\rho \leq \hat{\rho}_0$.} 
\begin{proof} Take  $x\in B_K(x_0)$ arbitrarily. Then  $d_K(x_0,x)<\frac{\delta_K(x_0)}{3K}$. 
Since $\varphi (d_K(x_0,\zeta ))<\alpha \delta_K(x_0)$, we have
\begin{align*}
\varphi (d_K(x,\zeta ))&<\varphi (K[\frac{\delta_K(x_0)}{3K}+d_K(x_0,\zeta )])\leq  
\varphi (K[\frac{d_K(x_0,\zeta )}{3K}+d_K(x_0,\zeta )])\\
&<\varphi ((K+\frac{1}{3})d_K(x_0,\zeta ))\leq \varphi (2^{K+\frac{1}{3}}d_K(x_0,\zeta ))\\
&< C_0^{K+\frac{1}{3}}\varphi (d_K(x_0,\zeta ))<C_0^{K+\frac{1}{3}}\alpha \delta_K(x_0)\leq \frac{C_0^{K+\frac{1}{3}}3\alpha K}{2}\delta_K(x),
\end{align*}
provided that $2^{K+\frac{1}{3}}d_K(x_0,\zeta )\leq r_0$. But this surely holds, since $x_0\in \Gamma_{\varphi ,\rho}(\zeta , \alpha )$ and $\rho \leq \hat{\rho}_0\leq \frac{1}{2^{\alpha +1}}\varphi \left(\frac{r_0}{2^{K+\frac{1}{3}}}\right)$. Hence $x\in \Gamma_{\varphi ,\rho' }(\zeta ,\alpha')$.
\end{proof}
\subsection{Lemma} \emph{ Let $C_1'=C_0^{\frac{3\alpha K}{2}C_0^{K+\frac{1}{3}}+1}$ and $\alpha'=\frac{3\alpha K}{2}C_0^{K+\frac{1}{3}}$. Then for all $x\in \Omega_{\rho'}$, where $\rho'=(K+\frac{1}{3})\rho$, one has ${\tilde{\Gamma}}_{\varphi}(x,\alpha')\subset B_K(x,C_1'\varphi^{-1}(\delta_K(x)))$, provided $0<\rho \leq \hat{\rho}_0$.}
\begin{proof} Suppose that ${\tilde{\Gamma}}_{\varphi}(x,\alpha')\ne \emptyset$ and take  $\xi \in {\tilde{\Gamma}}_{\varphi} (x,\alpha' )$ arbitrarily. But then $x\in \Gamma_{\varphi}(\xi ,\alpha')$, that is $\varphi (d_K(x,\xi ))<\alpha'\delta_K(x)$. Hence
\begin{equation*}d_K(x,\xi )<\varphi^{-1}(\alpha'\delta_K(x))<\varphi^{-1}(2^{\alpha'+1}\delta_K(x))\leq C_0^{\alpha'+1}\varphi^{-1}(\delta_K(x)),
\end{equation*}
provided that $2^{\alpha'+1}(K+\frac{1}{3})\hat{\rho}_0\leq r_0$. But this holds, since, by assumption, 
\begin{equation*}
\hat{\rho}_0\leq \frac{r_0}{2^{\frac{3\alpha K}{2}C_0^{K+\frac{1}{3}}+1}(K+\frac{1}{3})}.
\end{equation*}
Thus $\xi \in B_K(x,C_1'\varphi^{-1}(\delta_K(x)))$.
\end{proof}
\subsection{Lemma} \emph{Let $\zeta \in \partial \Omega$ and  $x\in \Gamma_{\varphi ,\rho}(\zeta , \alpha )$. Let $C_1'$ and $C_2'$  be as above. Then for $C_3'=K\left( 1+\frac{C_0^{\alpha +1}}{C_1'C_2'}\right)$ one has}  
\begin{equation*} B_K(x,C_1'C_2'\varphi^{-1}(\delta_K(x)))\subset B_K(\zeta ,C_1'C_2'C_3'\varphi^{-1}(\delta_K(x))),\end{equation*}
\emph{provided $0<\rho \leq \hat{\rho}_0$.}
\begin{proof} Take $z\in B_K(x,C_1'C_2'\varphi^{-1}(\delta_K(x)))$ arbitrarily. Then clearly
\begin{equation*}
d_K(z,\zeta )\leq K[d_K(x,z)+d_K(x,\zeta )]<K[C_1'C_2'\varphi^{-1}(\delta_K(x))+d_K(x,\zeta )].
\end{equation*}
Now $x\in \Gamma_{\varphi ,\rho}(\zeta ,\alpha )$, thus $\varphi (d_K(x,\zeta ))<\alpha \delta_K(x)$, and also 
\begin{equation*}
d_K(x,\zeta )<\varphi^{-1}(\alpha \delta_K(x))\leq \varphi^{-1}(2^{\alpha +1}\delta_K(x))\leq C_0^{\alpha +1}\varphi^{-1}(\delta_K(x)),
\end{equation*}
provided that $2^{\alpha +1}\hat{\rho}_0\leq r_0$, which again  holds, since, by assumption, $\hat{\rho}_0\leq \frac{r_0}{2^{\alpha +1}}$. Hence,
\begin{align*}
d_K(z,\zeta )&<K[C_1'C_2'\varphi^{-1}(\delta_K(x))+C_0^{\alpha +1}\varphi^{-1}(\delta_K(x))]\\
&<K(C_1'C_2'+C_0^{\alpha +1})\varphi^{-1}(\delta_K(x))=
C_1'C_2'K\left(1+\frac{C_0^{\alpha +1}}{C_1'C_2'}\right)\varphi^{-1}(\delta_K(x)),
\end{align*}
and so $z\in B_K(\zeta ,C_1'C_2'C_3'\varphi^{-1}(\delta_K(x)))$.
\end{proof}

Then our result, an extension to Theorem~1, pp.~131-132, of \cite{Ri09}:

\subsection{Theorem} {\emph{Let $X$ be a locally uniformly homogeneous space satisfying the condition (1). Suppose that $d_K:X\times X \rightarrow [0,+\infty )$ is separately continuous and a Borel function. Let 
$\varphi :[0,+\infty )\rightarrow [0,+\infty )$ be an admissible function, with constants $r_0$ and $C_0$. Let $\alpha >0$, $\gamma\in {\mathbb{R}}$, $d>0$ and $C_4>0$ be arbitrary. Let $u:X \rightarrow [0,+\infty )$ be a $K_1$-quasinearly subharmonic function. Then there is a constant $C=C(\alpha , \gamma , \varepsilon_0,d,A,C_0,C_4, K,K_1)$ such that for each component $X_1$ of $X$ and for each  domain $\Omega \subset X_1$, $\Omega \ne X_1$, whose   boundary $\partial \Omega$ is Ahlfors-regular from above, with dimension $d$ and with constant $C_4$,   
\begin{equation} \int\limits_{\partial \Omega }
\sup_{x\in \Gamma _{\varphi ,\rho }(\zeta ,\alpha )}\{\, \delta_K (x)^{\gamma }\mu (B_K(x))[\varphi ^{-1}(\delta_K (x))]^{-d}u(x)\,\}
\,d{\mathcal{H}}_K^d(\zeta )\leq C\int\limits_{\Omega_{\rho'}  }\delta_K (x)^\gamma\,u(x)\,d\mu (x),
\end{equation}
for all $\rho$, $0 < \rho \leq \hat{\rho}_0$.  
Here $\rho '=(K+\frac{1}{3})\rho $ and $\hat{\rho}_0$ is as above.} 
\subsubsection{Remark} Above and below we use the following (maybe unstandard, but in the considered situation of nonnegative functions nevertheless natural) convention: Let $A\subset X$, $B \subset A$ and  $g:A\rightarrow [0,+\infty ]$. If $B=\emptyset$,  then we \emph{define} $\sup_{x\in B}\{\, g(x)\,\}=0.$ 
\subsubsection{Remark} Though the constant 
$C$ above in (5) does depend on $K_1$ and on $\varepsilon_0$, it is, nevertheless, otherwise independent of the $K_1$-quasinearly subharmonic \mbox{function $u$.}
\begin{proof} Suppose $0<\rho \leq \hat{\rho}_0$. Write 
\begin{equation*}
E:=\{\, \zeta \in \partial \Omega \, :\, \Gamma _{\varphi ,\rho }(\zeta ,\alpha )\ne \emptyset \, \}.
\end{equation*}
Using the fact that $d_K(\cdot ,\cdot )$ is separately upper semicontinuous, one sees easily that $E$ is open in $\partial \Omega$.

Take $\zeta\in E$  and $x_0\in \Gamma_{\varphi ,\rho }(\zeta ,\alpha)$ arbitrarily.  Since $u$ is quasinearly subharmonic and 
\begin{equation*}
\frac{\epsilon_0\delta_K(x_0)}{3K}<\epsilon_0\delta_K(x_0)< \delta_K(x_0)\leq \rho \leq \hat{\rho}_0\leq \rho_0
\end{equation*}
one obtains
\begin{align*}u(x_0)&\leq \frac{K_1}{\mu (B_K(x_0,\frac{\epsilon_0\delta_K(x_0)}{3K}))}\int\limits_{B_K(x_0,\frac{\epsilon_0\delta_K(x_0)}{3K})}u(x)\,d\mu (x)\\
&\leq \frac{K_1}{\mu (B_K(x_0,\frac{\epsilon_0\delta_K(x_0)}{3K}))}\int\limits_{B_K(x_0,\frac{\delta_K(x_0)}{3K})}u(x)\,d\mu (x).
\end{align*}
Choose  $n_0\in {\mathbb{N}}$ such that 
\begin{equation*}2^{n_0-1}< \frac{1}{\epsilon_0}\leq 2^{n_0}.\end{equation*}
Then  
\begin{align*}\mu (B_K(x_0,\frac{\delta_K(x_0)}{3K}))&=\mu (B_K(x_0,\frac{1}{\epsilon_0}\cdot \frac{\epsilon_0\delta_K(x_0)}{3K}))
\leq \mu (B_K(x_0,2^{n_0}\cdot \frac{\epsilon_0\delta_K(x_0)}{3K}))\\
&\leq A^{n_0}\mu (B_K(x_0, \frac{\epsilon_0\delta_K(x_0)}{3K}))
\leq A^{1-\log_2\epsilon_0}\mu (B_K(x_0, \frac{\epsilon_0\delta_K(x_0)}{3K})).
\end{align*}
Hence 
\begin{equation}u(x_0)\leq \frac{K_1A^{1-log_2 \epsilon_0}}{\mu (B_K(x_0,\frac{\delta_K(x_0)}{3K}))}\int\limits_{B_K(x_0,\frac{\delta_K(x_0)}{3K})}u(x)\,d\mu (x)=\frac{K_1A^{1-log_2 \epsilon_0}}
{\mu (B_K(x_0))}\int\limits_{B_K(x_0)}u(x)\,d\mu (x).
\end{equation}
With the aid of the fact that
$\delta_K(x_0)^{\gamma}\leq (\frac{3K}{2})^{\mid \gamma \mid}\delta_K(x)^{\gamma}$ for all $x\in B_K(x_0)$ and that  $[\varphi ^{-1}(\delta_K(x_0))]^d\geq \frac{1}{C_0^{d(K+\frac{1}{3})}}[\varphi ^{-1}(\delta_K(x))]^d$ for all $x\in B_K(x_0)$,
we get, with the aid of  Lemma~2.1 above, from (6) above,
for all $C_1\geq 1$, 
\begin{align*}
&\frac{\delta_K(x_0)^{\gamma}\mu (B_K(x_0))u(x_0)}{[\varphi^{-1}(\delta_K(x_0))]^{d}+{\mathcal{H}}_K^d(B_K(x_0,C_1C_2'\varphi^{-1}(\delta_K(x_0)))\cap \partial  \Omega)}\leq\\
&\leq {K_1A^{1-log_2 \epsilon_0}}\int\limits_{B_K(x_0)}
\frac{\delta_K(x_0)^{\gamma}u(x)}{[\varphi^{-1}(\delta_K(x_0))]^{d}+{\mathcal{H}}_K^d(B_K(x_0,C_1C_2'\varphi^{-1}(\delta_K(x_0)))\cap \partial  \Omega)}\,d\mu (x)\leq \\
&\leq {K_1A^{1-log_2 \epsilon_0}}\int\limits_{B_K(x_0)}
\frac{\left( \frac{3K}{2}\right)^{\mid \gamma \mid}\delta_K(x)^{\gamma}u(x)}{\frac{1}{C_0^{d(K+\frac{1}{3})}}[\varphi^{-1}(\delta_K(x))]^{d}+{\mathcal{H}}_K^d(B_K(x,C_1\varphi^{-1}(\delta_K(x)))\cap \partial  \Omega)}\,d\mu (x)\leq\\
&\leq {\left(\frac{3K}{2}\right)^{\mid \gamma \mid}K_1A^{1-log_2 \epsilon_0}C_0^{d(K+\frac{1}{3})}}\int\limits_{B_K(x_0)}
\frac{\delta_K(x)^{\gamma}u(x)}{[\varphi^{-1}(\delta_K(x))]^{d}+{\mathcal{H}}_K^d(B_K(x,C_1\varphi^{-1}(\delta_K(x)))\cap \partial  \Omega )}\,d\mu (x).
\end{align*}
By Lemma~2.2, $B_K(x_0)\subset \Gamma_{\varphi ,\rho'}(\zeta,\alpha' )$, where $\rho'=(K+\frac{1}{3})\rho$ and $\alpha'=\frac{3\alpha K}{2}C_0^{K+\frac{1}{3}}$. Thus
\begin{align*}&\frac{\delta_K(x_0)^{\gamma}\mu (B_K(x_0))u(x_0)}{[\varphi^{-1}(\delta_K(x_0))]^{d}+{\mathcal{H}}_K^d(B_K(x_0,C_1C_2'\varphi^{-1}(\delta_K(x_0)))\cap \partial  \Omega)}\leq\\
&\leq {\left(\frac{3K}{2}\right)^{\mid \gamma \mid}K_1A^{1-log_2 \epsilon_0}C_0^{d(K+\frac{1}{3})}}\int\limits_{\Gamma_{\varphi ,\rho'}(\zeta , \alpha' )}
\frac{\delta_K(x)^{\gamma}u(x)}{[\varphi^{-1}(\delta_K(x))]^{d}+{\mathcal{H}}_K^d(B_K(x,C_1\varphi^{-1}(\delta_K(x)))\cap \partial  \Omega )}\,d\mu (x).
\end{align*}
Taking then the supremum on the left hand side over $x_0\in \Gamma_{\varphi ,\rho }(\zeta ,\alpha )$, we get
\begin{align*}&\sup_{x_0\in \Gamma_{\varphi ,\rho}(\zeta ,\alpha )} \frac{\delta_K(x_0)^{\gamma}\mu (B_K(x_0))u(x_0)}{[\varphi^{-1}(\delta_K(x_0))]^{d}+{\mathcal{H}}_K^d(B_K(x_0,C_1C_2'\varphi^{-1}(\delta_K(x_0)))\cap \partial  \Omega)}\leq\\
&\leq {\left(\frac{3K}{2}\right)^{\mid \gamma \mid}K_1A^{1-log_2 \epsilon_0}C_0^{d(K+\frac{1}{3})}}\int\limits_{\Gamma_{\varphi ,\rho'}(\zeta , \alpha' )}
\frac{\delta_K(x)^{\gamma}u(x)}{[\varphi^{-1}(\delta_K(x))]^{d}+{\mathcal{H}}_K^d(B_K(x,C_1\varphi^{-1}(\delta_K(x)))\cap \partial  \Omega )}\,d\mu (x).
\end{align*}
Next  integrate on both sides with respect to $\zeta$ over $E$ and use Fubini's theorem:
\begin{align*}&\int\limits_{E}\sup_{x_0\in \Gamma_{\varphi ,\rho}(\zeta ,\alpha )} \frac{\delta_K(x_0)^{\gamma}\mu (B_K(x_0))u(x_0)}{[\varphi^{-1}(\delta_K(x_0))]^{d}+{\mathcal{H}}_K^d(B_K(x_0,C_1C_2'\varphi^{-1}(\delta_K(x_0)))\cap \partial  \Omega)}\,d{\mathcal{H}}_K^d(\zeta )\leq\\
&\leq {\left(\frac{3K}{2}\right)^{\mid \gamma \mid}K_1A^{1-log_2 \epsilon_0}C_0^{d(K+\frac{1}{3})}}\int\limits_{E}\{\int\limits_{\Gamma_{\varphi ,\rho'}(\zeta , \alpha' )}
\frac{\delta_K(x)^{\gamma}u(x)}{[\varphi^{-1}(\delta_K(x))]^{d}+{\mathcal{H}}_K^d(B_K(x,C_1\varphi^{-1}(\delta_K(x)))\cap \partial  \Omega )}\,d\mu (x)\}d{\mathcal{H}}_K^d(\zeta )\\
&\leq {\left(\frac{3K}{2}\right)^{\mid \gamma \mid}K_1A^{1-log_2 \epsilon_0}C_0^{d(K+\frac{1}{3})}}\int\limits_{E}\{\int\limits_{\Omega_{\rho'}}\chi_{\Gamma_{\varphi}(\zeta , \alpha' )}(x)
\frac{\delta_K(x)^{\gamma}u(x)}{[\varphi^{-1}(\delta_K(x))]^{d}+{\mathcal{H}}_K^d(B_K(x,C_1\varphi^{-1}(\delta_K(x)))\cap \partial  \Omega )}\,d\mu (x)\}d{\mathcal{H}}_K^d(\zeta )\\
&\leq {\left(\frac{3K}{2}\right)^{\mid \gamma \mid}K_1A^{1-log_2 \epsilon_0}C_0^{d(K+\frac{1}{3})}}\int\limits_{E}\{\int\limits_{\Omega_{\rho'}}\chi_{{\tilde{\Gamma}}_{\varphi}(x, \alpha' )}(\zeta )
\frac{\delta_K(x)^{\gamma}u(x)}{[\varphi^{-1}(\delta_K(x))]^{d}+{\mathcal{H}}_K^d(B_K(x,C_1\varphi^{-1}(\delta_K(x)))\cap \partial  \Omega )}\,d\mu (x)\}d{\mathcal{H}}_K^d(\zeta )\\
&\leq {\left(\frac{3K}{2}\right)^{\mid \gamma \mid}K_1A^{1-log_2 \epsilon_0}C_0^{d(K+\frac{1}{3})}}\int\limits_{\Omega_{\rho'}}\{\int\limits_{E}\chi_{{\tilde{\Gamma}}_{\varphi}(x, \alpha' )}(\zeta )d{\mathcal{H}}_K^d(\zeta )\}
\frac{\delta_K(x)^{\gamma}u(x)}{[\varphi^{-1}(\delta_K(x))]^{d}+{\mathcal{H}}_K^d(B_K(x,C_1\varphi^{-1}(\delta_K(x)))\cap \partial  \Omega )}\,d\mu (x)\\
&\leq {\left(\frac{3K}{2}\right)^{\mid \gamma \mid}K_1A^{1-log_2 \epsilon_0}C_0^{d(K+\frac{1}{3})}}\int\limits_{\Omega_{\rho'}}{\mathcal{H}}_K^d({{\tilde{\Gamma}}_{\varphi}(x, \alpha' )})
\frac{\delta_K(x)^{\gamma}u(x)}{[\varphi^{-1}(\delta_K(x))]^{d}+{\mathcal{H}}_K^d(B_K(x,C_1\varphi^{-1}(\delta_K(x)))\cap \partial  \Omega )}\,d\mu (x).
\end{align*}
Choosing $C_1=C_1'$ and using then Lemma~2.3 we get
\begin{align*}&\int\limits_{E}\sup_{x_0\in \Gamma_{\varphi ,\rho}(\zeta ,\alpha )} \frac{\delta_K(x_0)^{\gamma}\mu (B_K(x_0))u(x_0)}{[\varphi^{-1}(\delta_K(x_0))]^{d}+{\mathcal{H}}_K^d(B_K(x_0,C_1'C_2'\varphi^{-1}(\delta_K(x_0)))\cap \partial  \Omega)}\,d{\mathcal{H}}_K^d(\zeta )\leq\\
&\leq {\left(\frac{3K}{2}\right)^{\mid \gamma \mid}K_1A^{1-log_2 \epsilon_0}C_0^{d(K+\frac{1}{3})}}\int\limits_{\Omega_{\rho'}}\frac{{\mathcal{H}}_K^d(B_K(x,C_1'\varphi^{-1}(\delta_K(x)))\cap \partial \Omega) 
\delta_K(x)^{\gamma}u(x)}{[\varphi^{-1}(\delta_K(x))]^{d}+{\mathcal{H}}_K^d(B_K(x,C_1'\varphi^{-1}(\delta_K(x)))\cap \partial  \Omega )}\,d\mu (x)\\
&\leq {\left(\frac{3K}{2}\right)^{\mid \gamma \mid}K_1A^{1-log_2 \epsilon_0}C_0^{d(K+\frac{1}{3})}}\int\limits_{\Omega_{\rho'}}\delta_K(x)^{\gamma}u(x)\,d\mu (x).
\end{align*}
On the other hand, by Lemma~2.4,
we get
\begin{align*}
&\sup_{x\in \Gamma_{\varphi ,\rho}(\zeta ,\alpha )} \frac{\delta_K(x)^{\gamma}\mu (B_K(x))u(x)}{[\varphi^{-1}(\delta_K(x))]^{d}+{\mathcal{H}}_K^d(B_K(x,C_1'C_2'\varphi^{-1}(\delta_K(x)))\cap \partial  \Omega)}\geq\\
&\geq \sup_{x\in \Gamma_{\varphi ,\rho}(\zeta ,\alpha )} \frac{\delta_K(x)^{\gamma}\mu (B_K(x))u(x)}{[\varphi^{-1}(\delta_K(x))]^{d}+{\mathcal{H}}_K^d(B_K(\zeta ,C_1'C_2'C_3'\varphi^{-1}(\delta_K(x)))\cap \partial  \Omega)}.
\end{align*}
Since $\partial \Omega$ is Ahlfors-regular from above, one has
\begin{equation*}{\mathcal{H}}_K^d(B_K(\zeta ,C_1'C_2'C_3'\varphi^{-1}(\delta_K(x)))\cap \partial  \Omega )\leq 
C_4[C_1'C_2'C_3'\varphi^{-1}(\delta_K(x))]^d<+\infty .
\end{equation*}
Therefore
\begin{align*}
&\sup_{x\in \Gamma_{\varphi ,\rho}(\zeta ,\alpha )} \frac{\delta_K(x)^{\gamma}\mu (B_K(x))u(x)}{[\varphi^{-1}(\delta_K(x))]^{d}+{\mathcal{H}}_K^d(B_K(\zeta ,C_1'C_2'C_3'\varphi^{-1}(\delta_K(x)))\cap \partial  \Omega)}\geq\\
&\geq \sup_{x\in \Gamma_{\varphi ,\rho}(\zeta ,\alpha )} \frac{\delta_K(x)^{\gamma}\mu (B_K(x))u(x)}{[\varphi^{-1}(\delta_K(x))]^{d}+C_4[C_1'C_2'C_3'\varphi^{-1}(\delta_K(x))]^d}=\\
&= \frac{1}{1+C_4(C_1'C_2'C_3')^d}\sup_{x\in \Gamma_{\varphi ,\rho}(\zeta ,\alpha )} \{\delta_K(x)^{\gamma}\mu (B_K(x)) [\varphi^{-1}(\delta_K(x))]^{-d} u(x)\}.
\end{align*}
Thus we have:
\begin{equation}
\int\limits_{E} \sup_{x\in \Gamma_{\varphi ,\rho}(\zeta ,\alpha )} \{\delta_K(x)^{\gamma}\mu (B_K(x)) [\varphi^{-1}(\delta_K(x))]^{-d} u(x)\}d{\mathcal{H}}_K^d(\zeta )\leq C\,\int\limits_{\Omega_{\rho'}}\delta_K(x)^{\gamma}u(x)\, d\mu (x),
\end{equation}
where
\begin{equation*}
C=\left(\frac{3K}{2}\right)^{\mid \gamma \mid}K_1A^{1-log_2 \epsilon_0}C_0^{d(K+\frac{1}{3})}[1+C_4\left(C_1'C_2'C_3'\right)^d]
\end{equation*}
and
\begin{equation*}
C_1'=C_0^{\frac{3\alpha K}{2}C_0^{K+\frac{1}{3}}+1}, \quad C_2'=\frac{C_0^{\alpha}}{3}+KC_0^{K+\frac{1}{3}},\quad C_3'=K\left( 1+\frac{C_0^{\alpha +1}}{C_1'C_2'}\right).
\end{equation*}

To conclude the proof,  observe  the following. 
First, since $\Gamma_{\varphi ,\rho}(\zeta ,\alpha )=\emptyset$ for all $\zeta \in \partial \Omega\setminus E$, we can, just using our convention in Remark~2.5.1, replace (7) by the desired inequality:
\begin{equation*}
\int\limits_{\partial \Omega} \sup_{x\in \Gamma_{\varphi ,\rho}(\zeta ,\alpha )} \{\delta_K(x)^{\gamma}\mu (B_K(x)) [\varphi^{-1}(\delta_K(x))]^{-d} u(x)\}d{\mathcal{H}}_K^d(\zeta )\leq C\,\int\limits_{\Omega_{\rho'}}\delta_K(x)^{\gamma}u(x)\, d\mu (x).
\end{equation*}

Second, the functions  
\begin{equation*}
\partial \Omega \ni \zeta \mapsto \sup_{x\in \Gamma_{\varphi ,\rho}(\zeta ,\alpha )} \frac{\delta_K(x)^{\gamma}\mu (B_K(x))u(x)}{[\varphi^{-1}(\delta_K(x))]^{d}+{\mathcal{H}}_K^d(B_K(x,C_1C_2'\varphi^{-1}(\delta_K(x)))\cap \partial  \Omega)}\in [0,+\infty ]
\end{equation*}
and
\begin{equation*}
\partial \Omega \ni \zeta \mapsto \sup_{x\in \Gamma_{\varphi ,\rho}(\zeta ,\alpha )} \{\delta_K(x)^{\gamma}\mu (B_K(x)) [\varphi^{-1}(\delta_K(x))]^{-d} u(x)\}\in [0,+\infty ]
\end{equation*}
are lower semicontinuous. Thus the above integrations on ``the left hand sides'' are justified. 

Third, the functions 
\begin{align*}
\Omega_{\rho'} \ni x&\mapsto  
\frac{\delta_K(x)^{\gamma}u(x)}{[\varphi^{-1}(\delta_K(x))]^{d}+{\mathcal{H}}_K^d(B_K(x,C_1\varphi^{-1}(\delta_K(x)))\cap \partial  \Omega )}\in [0,+\infty )
\end{align*}
and 
\begin{align*}
\Omega_{\rho'}\times \partial \Omega \ni (x,\zeta )&\mapsto  \chi_{{{\Gamma}}_{\varphi}(\zeta , \alpha' )}(x)
\frac{\delta_K(x)^{\gamma}u(x)}{[\varphi^{-1}(\delta_K(x))]^{d}+{\mathcal{H}}_K^d(B_K(x,C_1\varphi^{-1}(\delta_K(x)))\cap \partial  \Omega )}=\\
&=  \chi_{{\tilde{\Gamma}}_{\varphi}(x, \alpha' )}(\zeta )
\frac{\delta_K(x)^{\gamma}u(x)}{[\varphi^{-1}(\delta_K(x))]^{d}+{\mathcal{H}}_K^d(B_K(x,C_1\varphi^{-1}(\delta_K(x)))\cap \partial  \Omega )}\in [0,+\infty )
\end{align*}
are Borel measurable. Hence the integrations  and the use of Fubini's theorem on ``the right hand sides'' are justified, too. Observe that here we  use  our additional assumption that the K-quasimetric 
$d_K$ is  Borel measurable.
\end{proof}
\subsubsection{Remark} At present we do not know whether our assumption that the $K$-quasimetric $d_K:X\times X \rightarrow [0,+\infty )$ is separately continuous and Borel measurable, is really necessary or not. Observe anyway that a quasimetric $d_K$ is separately upper semicontinuous, see (iii) 4$^o$ above. But a separately upper semicontinuous function need not, however,  be measurable, see \cite{Si20} and  e.g. \cite{Gr02}, Example~1,  p.~11.  On the other hand,  if $K=1$, that is, if $X$ is a metric space, and the  function $d_1$ is separately continuous, then $d_1$ is Borel measurable by a result of Kuratowski, see \cite{Ru81}, p.~742, and the references therein, say.  Observe also that if the considered locally uniformly homogeneous space $X$ is moreover locally compact, then by a result of Johnson, see \cite{Jo69}, Theorem~2.2, p.~422, and again \cite{Ru81}, p.~742, $d_K$ is indeed Borel measurable.
\subsection{Corollary} {\emph{Let $X$ be a locally uniformly homogeneous space  satisfying the condition (1). Suppose that $d_K:X\times X \rightarrow [0,+\infty )$ is separately continuous and a Borel function. Let 
$\varphi :[0,+\infty )\rightarrow [0,+\infty )$ be an admissible function, with constants $r_0$ and $C_0$. Let $\alpha >0$, $\gamma \in {\mathbb {R}}$, $d>0$ and $C_4>0$ be arbitrary. Let $u:X \rightarrow [0,+\infty )$ be a $K_1$-quasinearly subharmonic function. Suppose that  $\Omega \subset X$, $\Omega \ne X$, is a domain whose    boundary $\partial \Omega$ is Ahlfors-regular from above, with dimension $d$ and with constant $C_4$, and that   
\begin{equation} \int\limits_{\Omega}\delta_K (x)^\gamma\,u(x)\,d\mu (x)<+\infty .
\end{equation}
Then for ${\mathcal{H}}_K^d$-almost every $(\varphi ,\alpha )$-accessible point $\zeta \in \partial \Omega$,
\begin{equation*}
\lim_{\rho \rightarrow 0}\left(\sup_{x\in \Gamma _{\varphi ,\rho}(\zeta ,\alpha )}\{\, \delta_K(x)^{\gamma }\mu (B_K(x))
[\varphi ^{-1}(\delta_K(x))]^{-d}u(x)\,\}\right)=0.
\end{equation*}}}
\subsubsection{Remark} If, instead of a locally uniformly homogeneous space $X$, one works in an Euclidean space ${\mathbb{R}}^n$, $n\geq 2$,  then   slightly better results hold: Namely,  one can omit the assumed Ahlfors-regularity condition of $\partial \Omega$, see the already cited results  \cite{Ri00}, Theorem, p.~233,  \cite{Ri03}, Theorem~2, pp.~175--176, \cite{Ri04}, Theorem~3.4.1, pp.~198--199, \cite{Ri06$_1$}, Theorem, p.~31, and \cite{PaRi08}, Theorem~4, p.~102. Observe, however, that the possibility for this omission in the Euclidean setup is based essentially on a well-known density estimate result for Hausdorff measures (which in turn is based, among others, on Vitali's covering theorem and  thus is of ``very Euclidean space-type''), see e.g. \cite{Ma95}, Theorem~6.2, p.~89. 
\begin{proof} By Theorem~2.5,
\begin{equation*} \int\limits_{\partial \Omega }
\sup_{x\in \Gamma _{\varphi ,\rho }(\zeta ,\alpha )}\{\, \delta_K (x)^{\gamma }\mu (B_K(x))[\varphi ^{-1}(\delta_K (x))]^{-d}u(x)\,\}
\,d{\mathcal{H}}_K^d(\zeta )\leq C\int\limits_{\Omega_{\rho'}  }\delta_K (x)^\gamma\,u(x)\,d\mu (x),
\end{equation*}
where $C=C(\alpha , \gamma , \varepsilon_0,d,A,C_0,C_4, K,K_1)$. Using then just Fatou's lemma and (8), one sees that 
\begin{equation*} \int\limits_{\partial \Omega }
\liminf_{\rho \rightarrow  0} \left(\sup_{x\in \Gamma _{\varphi ,\rho }(\zeta ,\alpha )}\{\, \delta_K (x)^{\gamma }\mu (B_K(x))[\varphi ^{-1}(\delta_K (x))]^{-d}u(x)\,\}\right)
\,d{\mathcal{H}}_K^d(\zeta )\leq C\,\liminf_{\rho \rightarrow 0}\int\limits_{\Omega_{\rho'}  }\delta_K (x)^\gamma\,u(x)\,d\mu (x) =0.
\end{equation*}
Thus the claim follows.
\end{proof}

\end{document}